\documentclass[prl,twoside,onecolumn]{revtex4}
\usepackage[utf8]{inputenc}
\usepackage{graphicx}

\date{\today}

\begin{document}
\title{The tightest knot is not necessarily the smallest}
\author{Alexander R. Klotz }
\affiliation{Department of Physics and Astronomy, California State University, Long Beach}
\begin{abstract}
In this note, we attempt to find counterexamples to the conjecture that the ideal form of a knot, that which minimizes its contour length while respecting a no-overlap constraint, also minimizes the volume of the knot, as determined by its convex hull. We measure the convex hull volume of knots during the length annealing process, identifying local minima in the hull volume that arise due to buckling and symmetry breaking. We use T(p,2) torus knots as an illustrative example of a family of knots whose locally minimal-length embeddings are not necessarily ordered by volume. We identify several knots whose central curve has a convex hull volume that is not minimized in the ideal configuration, and find that $8_{19}$ has a non-ideal global minimum in its convex hull volume even when the thickness of its tube is taken into account. 
    
\end{abstract}

\maketitle

\section{Introduction}

Physical knot theory is concerned with the properties of knots and links that are subject to physical constraints such as finite thickness or bending stiffness. One of the most studied quantities in physical knot theory is the ropelength of a knot, which is the minimum contour length of a knot represented by a tube of unit radius that cannot overlap with itself. The embedding of a knot that minimizes contour length is known as \textit{ideal} \cite{stasiak1998ideal}. Several theorems have been proven bounding the relationship between the crossing number of a knot and its ropelength \cite{diao2003lower, diao2020ropelength,buck1998four}, and computational studies suggest stronger bounds than those currently proven \cite{ashton2011knot,klotz, diao2006numerical}.

Many of the invariants found in physical knot theory are in some way related to the length of the knot: ropelength, ribbon length, lattice number, stick number. Algorithms have been developed and optimized to find or minimize these parameters. Since knots exist in three dimensions, it is natural to ask what the minimum volume of a knot is. In a recent work \cite{klotz}, we examined the ropelength of complex torus and satellite knots and posed the question of whether the tightest knot, the form with the smallest length-to-radius ratio, is also the smallest in terms of volume. In this work, we answer this question in the negative through illustrative counterexamples.

In addition to filaments and magnetic flux tubes, much of physical knot theory is motivated by knots that occur in DNA. Extremely tight knots are found in viral capsids \cite{arsuaga2008dna}, and it ultimately the volume of the knot that can determine whether it will fit within the capsid, not its ropelength. The diffusion of a knotted DNA molecule confined in a nanochannel is sped up by the effective shortening of the molecule, which depends on contour length, but slowed down by friction with the walls, which depends on the physical dimensions of the knot with respect to its confining geometry (experiments indicate that the shortening effect is stronger \cite{ma2021diffusion}). Knots have been detected in nanopore translocation experiments \cite{plesa, rajesh}, but as nanopore technology becomes more precise, it will be the cross-sectional area of the knot that affects the translocation signal, not its contour length.

Several studies have examined the volumes of tight knots. Millet and Rawdon examined the convex hull volume of ideal knots when normalized by contour length \cite{millett2003energy}, and Janse van Rensburg and Rechnitzer examined the bounding prisms and polyhedra of minimal lattice knots in order to calculate their compressiblity \cite{rensburg}. Celoria and Mahler recently developed a mapping between the topological spaces surrounding knot point sets and their physical volume, which they used to quantify a departure from ideality of a given knot embedding \cite{celoria}. To our knowledge, no algorithm has been developed to find the minimum-volume knot, and the above analyses focus on the volumes of minimum-length knots.

It is straightforward to demonstrate that a ropelength-minimizing configuration does not necessarily minimize volume. Consider a $2_{1}^{2}$\#$2_{1}^{2}$ Hopf chain (three links connected linearly), in which the first and third components are circles and the central component is a stadium curve (Figure 1). This is a ropelength minimizing configuration, but the two outer circles may be rotated to generate configurations with different volume at the same contour length. It is less straightforward to show that there are non-ideal configurations that have less volume than an ideal configuration.

\begin{figure}
    \centering
    \includegraphics[width=0.5\textwidth]{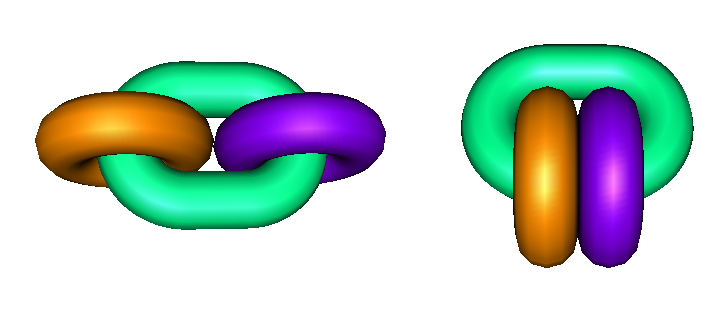}
    \caption{Two ideal configurations of a three-component Hopf chain. Both have a ropelength of $12\pi+4$, but the elongated configuration on the left occupies a rectangular prism of volume $6\times6\times12=432$ and the compact configuration on the right occupies a rectangular prism of volume $8\times8\times6=384$.}
    \label{fig:hopf}
\end{figure}

The volume taken up by a knot is ambiguous. The volume displaced by a physical knot submerged in liquid is simply its cross-sectional area ($\pi$, with unit radius) times the contour length, but this is not the most useful quantity if discussing, for example, the packing of knotted tubes within a finite volume. We may also consider the minimum-volume rectangle prism or ellipsoid that bounds a knot. In this work, we consider the volume of the convex hull of the knot as the most general and useful representation of the space taken up by a knot. There are two forms of the convex hull volume that we consider: the convex hull volume of the bare central curve of the knot, calculated based on its set of vertices, as well as the convex hull volume of the thick tube of the knot. The volume based on the vertices of the knots, which we refer to as the point hull volume, was considered by Millet and Rawdon \cite{millett2003energy}, who examined the hull volume of ideal knots when normalized by contour length.

Estimating the volume of the knot when its tube thickness is taken into account, which we refer to as the tube hull volume, may be done by calculating the parallel surface located a unit distance from the hull of the knot's central curve. While such an algorithm is not readily available, we measured the ``thick'' tube hull volume using a brute-force algorithm that generated a ring of 20 points in a circle around each vertex of the knot, normal to the tangent vector along the knot, and calculated the convex hull volume of the set of rings.

\section{Methodology}

Knots were annealed using the Constrained Gradient Optimization algorithm implemented in Ridgerunner, a software package developed and maintained by Jason Cantarella, and run on a virtual machine hosted by Amazon Web Services. A detailed description of Ridgerunnner's algorithm may be found in Ashton et al \cite{ashton2011knot}. Unlike our previous work \cite{klotz}, we annealed the knots without an equilateralization force, which did not result in pathological configurations due to the comparative simplicity of our knots and smooth initial conditions.

Our T(p,2) torus knots were initialized with 400 vertices, a minor radius of 2 and a major radius of $4+\frac{p}{\pi}$, which allows Ridgerunner to rescale the knot to a fairly tight configuration before annealing. Other torus knots ($8_{19}$ and $10_{124}$) were initialized with 200 vertices and the same toric radii, while non-torus knots were annealed from the initial configurations provided by Ridgerunner, which typically have less than 100 vertices.

To calculate convex hull volumes, we use the ``convhull'' function in MATLAB, which implements the quickhull algorithm developed by Barber et al. \cite{barber}. To calculate the tube hull volume, we encircle each knot vertex by 20 points in a plane perpendicular to the tangent vector between that vertex and its neighbor. The plane is defined by two orthogonal unit vectors. The first planar vector is the normalized cross product of the tangent vector and the radial vector from the knot's center of mass, and the second is the normalized cross product of the tangent vector and the first planar vector. The convex hull of the extended point cloud was then calculated using the same quickhull algorithm. If the number of points in each ring is increased from 20 to 40, the convex hull volume increases by approximately a third of a percent. The convex hulls are visualized with MATLAB's ``trisurf'' function, the produces triangulated surfaces shaded with a color scheme with the purpose of making surfaces easier to differentiate by eye. Knots are visualized using KnotPlot \cite{scharein2002interactive}. Figure 2 shows a tight configuration of a trefoil and its generated thick tube, along with the convex hull of the central curve (the point hull) and of the thick tube (the tube hull).

\begin{figure}
    \centering
    \includegraphics[width=0.5\textwidth]{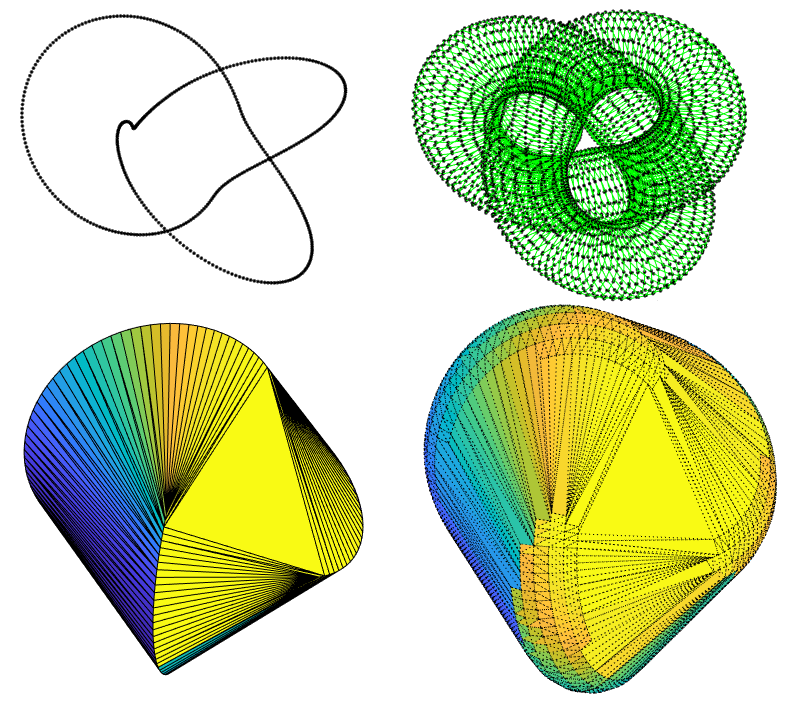}
    \caption{Tight configuration of a trefoil knot visualized as the points of its central curve (left) and its thick tube (right) generated from 20 points around each point on the central curve. Bottom images show the convex hulls of each set of points: we refer the hull of the central curve as the point hull and the hull of the tube as the tube hull.}
    \label{fig:tref}
\end{figure}

\section{Results and Discussion}

The first family of knots we examine are T(p,2) torus knots (for odd p), which are typically parameterized in the form of a circular double helix. We examine this family because it is known to have several families of locally minimal ropelength configurations and an annealing process that is characterized by spontaneous symmetry breaking. The contour length of the double helical configuration may be minimized with respect to its major and minor radii, into a locally tight form determined by Olsen and Bohr \cite{olsen2012geometry}. Klotz and Maldonado \cite{klotz} showed that for similar T(2k,2) torus links, the contour length when one of the components is a circle is tighter than a symmetric double helix configuration, with the circular component adopting a small helical pitch upon annealing. Earlier, Pieranski \cite{pieranski1998search} showed that a configuration tighter than the double helix exists, in which a helix is wound around an elongated tube, and the entire ``wrapped tube'' configuration is itself helically twisted. Recently, Huh et al. \cite{huh2018ropelength} showed that a superhelical configuration of T(6k+1,2) knots can be minimized with the respect to the parameters of the superhelix to find a locally minimal form that is tighter than the symmetric double helix for k$>2$, and more recently extended their analysis to 2-bridge knots \cite{huh2021tight}. Since the major radius of the torus must increase with p, the convex hull volume of the double helical configurations of T(p,2) knots grows quadratically with p, while the wrapped tube and twisted superhelix configurations grow linearly in p. We may hypothesize then that there are knots of intermediate complexity for which one minimized configuration has a lower contour length, while another has a lower convex hull volume. 



Initializing T(p,2) torus knots with their double helical parameterization, their contour steadily decreases but the symmetry of the double helix is maintained. When a critically tight double helix is reached, symmetry is spontaneously broken and the knot undergoes a transition analogous to a buckling instability. This is characterized by an elongation and contraction of the two axes in the initial plane of the double helix, and an expansion into the normal direction. A striking example may be seen on the late Piotr Pieranski's youtube channel \cite{youtube} and in the attached supplementary videos \footnote{The videos were generated using KnotPlot, not Ridgerunner, but the same features are observed. The ``tangential force'' used to reduce contour length in KnotPlot induces the rotation seen in the videos.}. While the contour length steadily decreases during the annealing process, the convex hull volume reaches a local minimum as symmetry is broken, before rising towards a local maximum. This non-monotonic hull volume with respect to contour length, may be seen for T(7,2) and T(13,2) in figures 3 and 4, but was observed in all T(p,2) knots we annealed for p$\geq7$. Eventually, a new minimum volume is reached in the ideal configuration, which for large p is a global minimum. However, the point hull volume of the T(7,2) knot (Fig. 3) is an exception: its helical minimum before symmetry breaking is the sole counterexample in this family. For larger p, the wrapped tube configuration has a smaller hull volume than the double helix, and the $3_1$ and $5_1$ knots are not sufficiently complex to undergo the same symmetry breaking. The volume of the locally minimal configurations is smaller when the knot contains more vertices, but it is unlikely that additional counterexamples may be found this way: with 200 vertices, the $9_1$ knot reaches a locally minimal point hull volume of 161.3, and with 600, it reaches 159.5. The global minimum is approximately 154.1.

Figure 3a shows the point and tube hull volumes as a function of the contour length at each stage of the annealing. In both cases, the helical minimum appears as a sharp cusp. In reality, the cusp can involve thousands of iterations, as the contour anneals through a very strong set of local minima before proceeding in its path towards the global minimum. Figure 3b shows the conformations at the helical volume minimum and at the ropelength minimum. The minimal-ropelength structure begins to resemble the wrapped tubes presented in previous works, but is not large enough to elongate and twist. When viewed from the side, it can be seen that the broken symmetry leads to a greater span in the direction normal to the initial circle of the torus. The convex hull of the helical configuration retains its heptagonal symmetry, while the convex hull of the ideal configuration resembles a Steinmetz solid or a g\"{o}mb\"{o}c.




\begin{figure}
    \centering
    \includegraphics[width=0.8\textwidth]{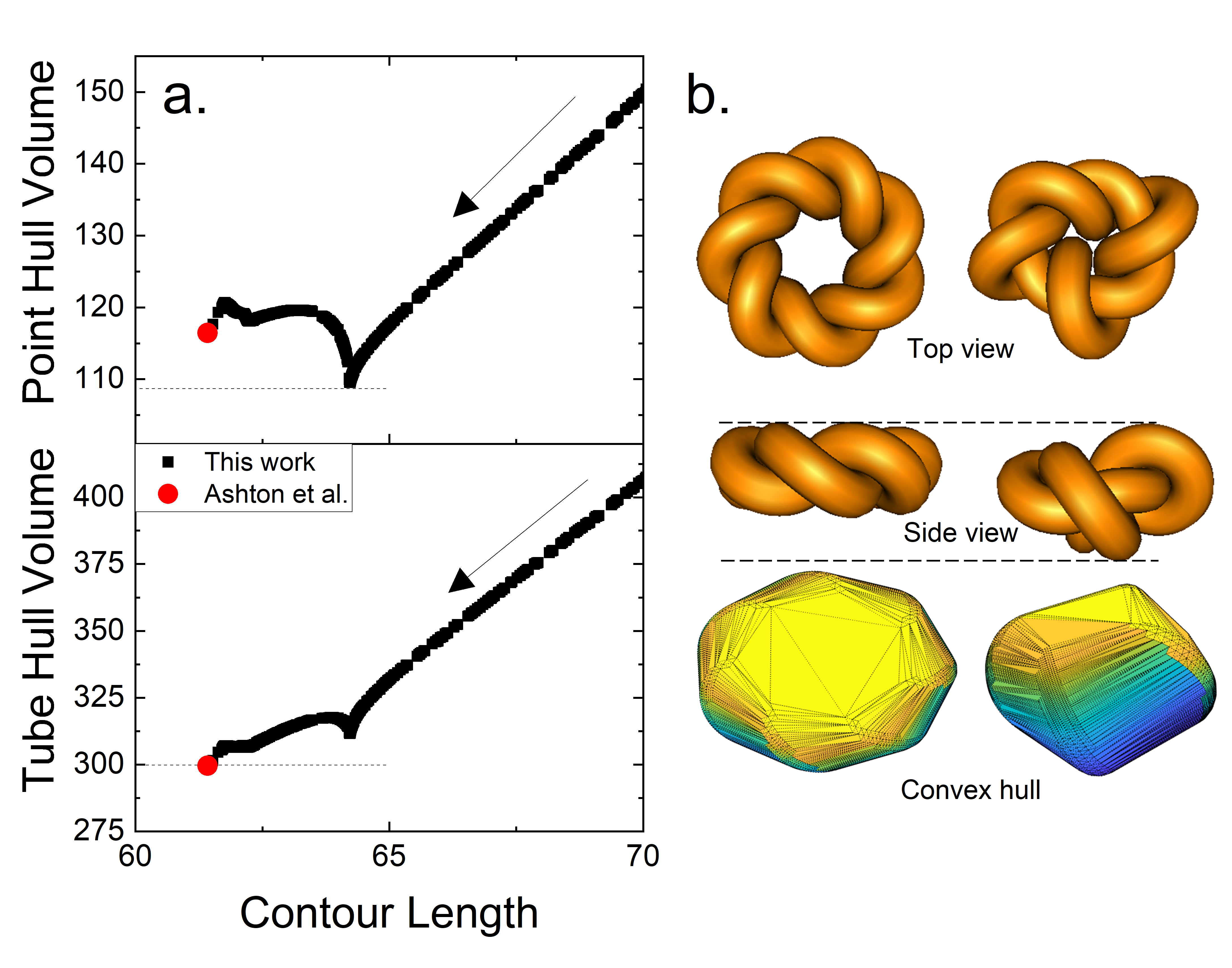}
    \caption{a. Convex hull volume of the $7_1$ knot as a function of its contour length during the annealing process, which proceeds from top-right to bottom-left on these charts. The top chart shows the convex hull volume of the central curve and the bottom figure shows the convex hull volume of the thick tube of the knot. Data from the ideal configurations of Ashton et al. \cite{ashton2011knot} are shown for comparison. In the upper chart, the global minimum in hull volume does not occur when contour length is minimized, but it does for the bottom chart. b. Visualizations of the point hull volume minimizing (left) and ropelength minimizing (right) configurations as well as their convex hulls.}
    \label{fig:71}
\end{figure}

The T(13,2) knot also presents an interesting case. Huh et al \cite{huh2018ropelength} proved that the minimal superhelix configuration of T(6k+1,2) knots has a smaller contour length than the minimal double helix of Olsen and Bohr for k$\geq$3  \cite{olsen2012geometry}. The superhelical T(7,2) and T(13,2) knots have a larger contour length than their helical counterparts due to the closures at the top and bottom of the helices. However, the superhelical configuration of T(13,2) has a significantly smaller convex hull volume, as measured from both the central curve and the thick tube, compared to the double helix (Figure 4a). The superhelix, the helical minimum, and the ideal wrapped tube all have approximately the same point hull volume: in Fig. 4a the superhelix has a point hull volume of 288.78, the minimal double helix has 289.27, and the wrapped tube 288.53. The superhelical configuration is not reached through annealing from a double helix initial condition, and highlights the sensitivity to initial conditions of these data. We also note that the Olsen and Bohr helical minimum lies on the annealing curve of our data, but does not coincide with the local minimum. This may imply an onset of broken symmetry before the dramatic buckling transition. The ideal configurations and their hulls are presented in Fig. 4b.


\begin{figure}
    \centering
    \includegraphics[width=0.8\textwidth]{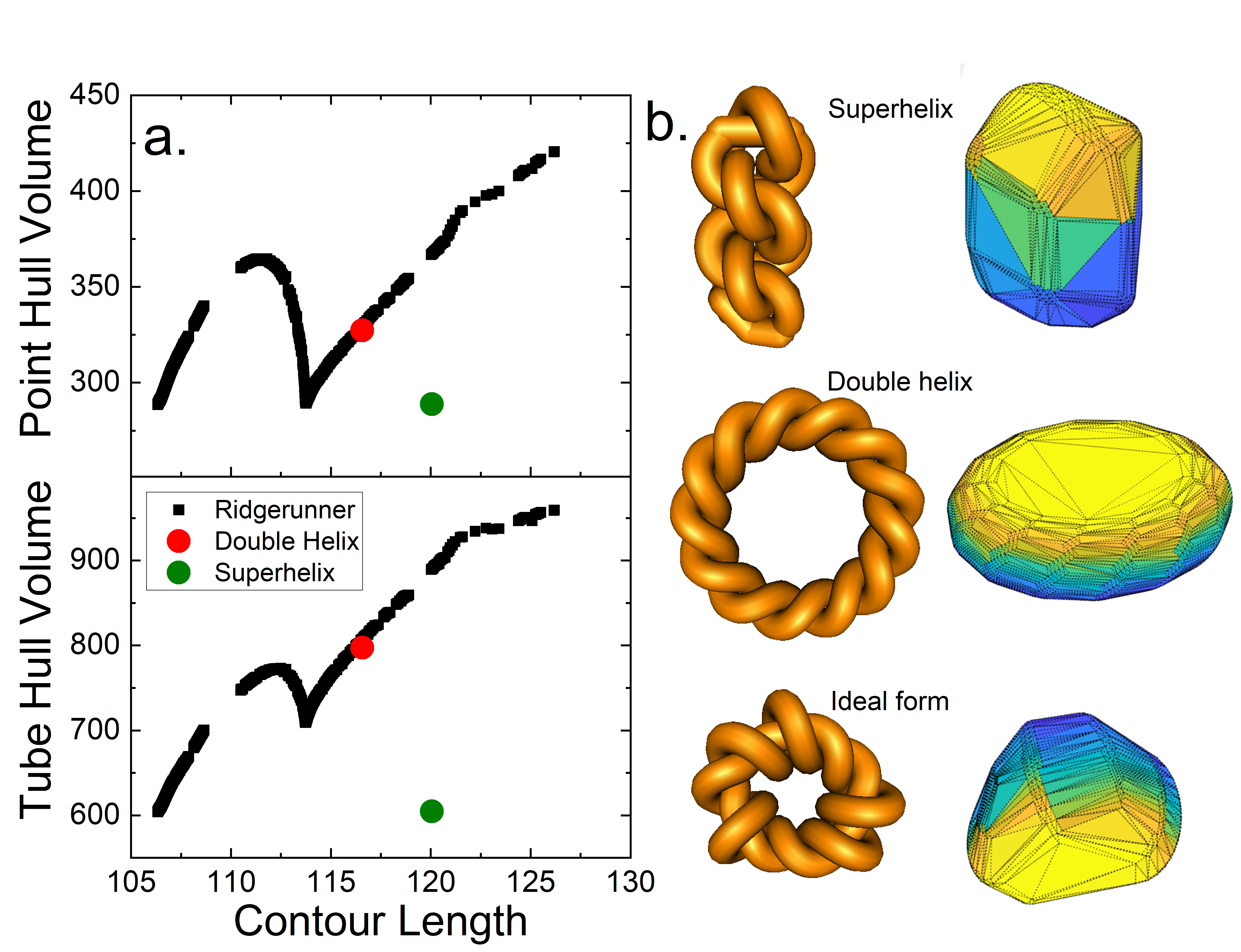}
    \caption{a. Convex hull volume of the T(13,2) knot as a function of its contour length during the annealing process, following the conventions of Fig. 3. The minimal double helix configuration from Olsen and Bohr \cite{olsen2012geometry} and the minimal superhelical configuration from Huh et al. \cite{huh2018ropelength} are shown. This knot is unique in that the superhelical form has a longer contour length but a smaller hull volume than the double helix, although neither form is globally minimal.  b. Visualizations of the minimal configurations configurations and their convex hulls. From top to bottom: superhelix, double helix, ideal wrapped tube.}
    \label{fig:13}
\end{figure}

Having found a counterexample in the T(p,2) family for the point hull volume, we annealed every knot with 3-7 crossings using Ridgerunner's default coordinates, and the first two non-alternating torus knots, $8_{19}$ and $10_{124}$, with a harmonic parameterization. The following knots were found to have minimal point hull volumes separate from their ideal configuration: $5_2$, $6_1$, $6_2$, $7_1$, $7_4$, $8_{19}$, and $10_{124}$. The only definitive non-ideal minimizer of tube hull volume was the $8_{19}$ knot (Fig. 5). Generally speaking, the minimal hull volumes decrease when more vertices are included in the knot, and it is possible that if more vertices were used in the annealings, more such counterexamples would be found. Many of these knots have an intricate pattern of minima and maxima in contour-volume space, which may warrant future study.

The trefoil's hull volume monotonically decreased through the annealing process, but curiously we observed the minimal tube hull volume of the trefoil knot to be below that of the ideal configuration of Ashton et al. \cite{ashton2011knot}, which was found with 2400 vertices compared to our 512, as well as Gilbert's ideal configuration \cite{katlas} found using the SONO (shrink-on-no-overlap) algorithm with 512 vertices. This may imply a globally minimal volume configuration with a contour length very close to its ropelength. The minimum contour length we reached in our annealing was 32.79, at which the hull volumes were 36.82 and 137.83, compared to 32.74, 37.51, and 138.88 from Ashton et al. \cite{ashton2011knot}. The tightest trefoil knot has been annealed with over 200,000 nodes \cite{przybyl2014high}; the raw data from that annealing may provide further insight.


As stated, the only knot for which the global minimum of the tube convex hull volume was found to be non-ideal was $8_{19}$, also known as T(4,3). Its annealing data may be seen in Figure 5. Along with the $10_{124}$/T(5,3), the ideal configuration appeared to have a volume which did not occupy even a local minimum. The tightest and smallest configurations are shown, with the only visual difference being the shifting of part of the tube's orientation at the top-left of each image.

\begin{figure}
    \centering
    \includegraphics[width=0.5\textwidth]{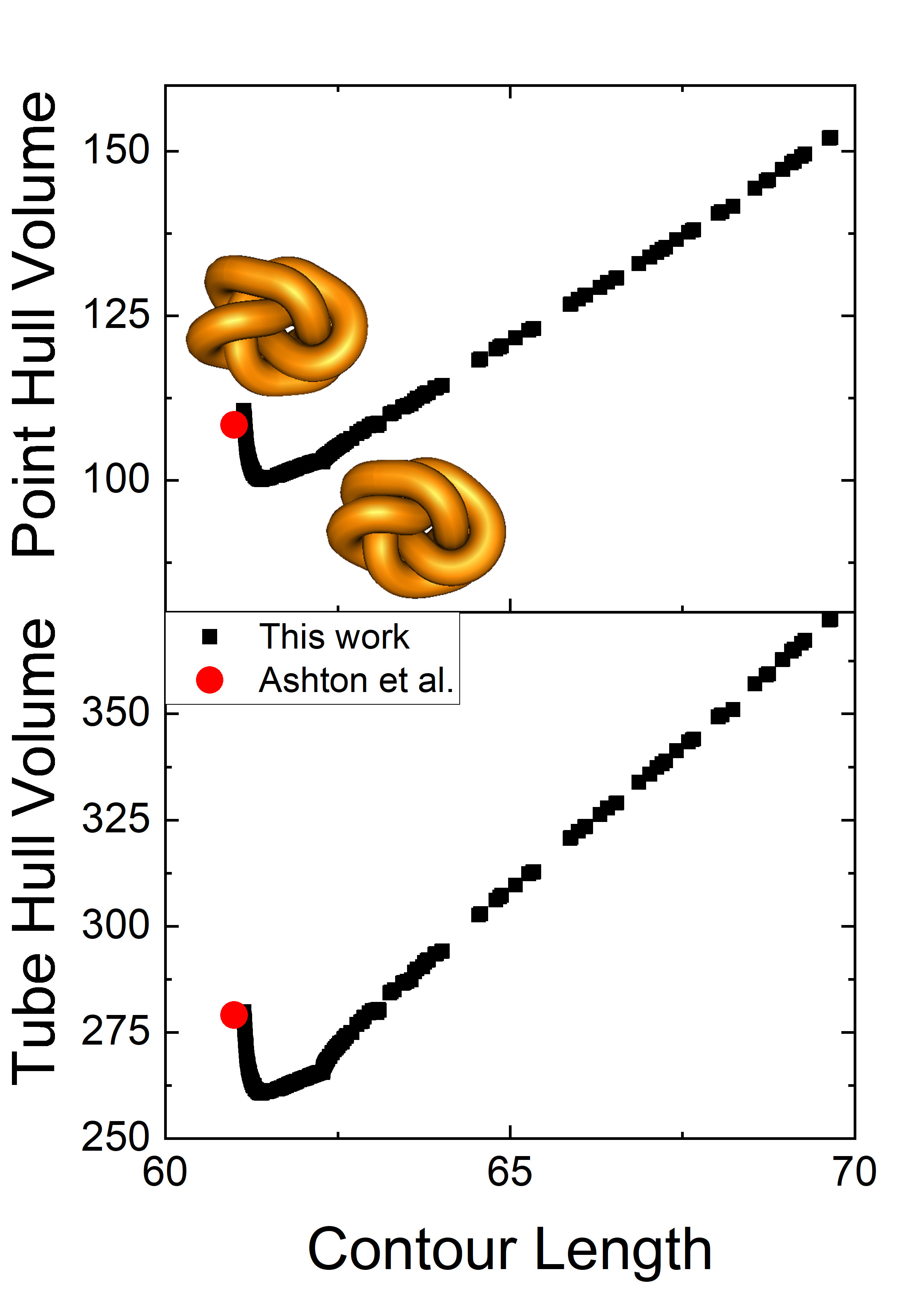}
    \caption{a. Convex hull volume of the $8_{19}$ knot as a function of its contour length during the annealing process, following the conventions of Fig. 3. Data from the ideal configurations of Ashton et al. \cite{ashton2011knot} are shown for comparison. Unlike the previous examples, the tube hull volume has a non-ideal global minimum. The ropelength-  and volume-minimizing configurations are shown, and only have minor visual differences.}
    \label{fig:819}
\end{figure}

A counterexample has been identified, but we must emphasize the assumption that the final configuration reached by ropelength annealing algorithms, whether in this work or those found by Ashton et al. using Ridgerunner \cite{ashton2011knot} and Gilbert using SONO \cite{katlas}, are close enough to the ``true'' ideal configuration such that our counterexamples are meaningful. This assumption may not be relaxed until proven lower bounds on ropelength are significantly strengthened: the proven lower bound for the length of the $7_1$ knot is 29.3 \cite{diao2003lower}, while the best measured upper bound is 61.4 \cite{ashton2011knot}.


Other physical parameters display interesting behavior during the ropelength annealing process. One such parameter is the volume fraction of the knot within its hull, $\pi\frac{L}{V}$. Empirically, this seems to plateau around 0.75, also close to the maximum packing density of spheres. Although this manuscript is concerned with volume, we may also examine the minimum cross-sectional area of a knot, which we can estimate as the geometric mean of the two smallest eigenvalues of the gyration tensor of the knot. This is a parameter that may be relevant to nanopore analyses of knotted DNA, in which the measured ionic current is displaced in proportion to the cross-sectional area of the object translocating through the pore. When non-ideal minima in the area exist (as for example in $8_{19}$), they tend to be stronger than those observed in the hull volume.

It is apparent from the data that the ratio of the tube hull volume to the point hull volume is not a constant. The ratio (for ideal knots) is somewhat analogous to a surface-to-volume ratio and generally decreases with increasing knot complexity. For the trefoil it reaches around 3.75, for the $10_{124}$ around 2.4, and for the T(19,2) around 1.9.

\section{Conclusion}

We sought to investigate the convex hull volume of knots during the ropelength minimization process in order to find counterexamples to the conjecture that the tightest configuration of a knot is also the smallest. A definitive counterexample ($8_{19}$) was identified, and in that sense our study was successful. Studying the broken symmetry of the double helical configuration of T(p,2) knots was illustrative in demonstrating the existence of local minima in hull volume. Comparison between the double helix and superhelix configurations showed that embeddings that are minimal with respect to the free parameters of their parameterization may have disordered length-volume relationships. Several partial counterexamples, minimizing the convex hull volume of the central curve but not the thick tube, were identified: $5_2$, $6_1$, $6_2$, $7_1$, $7_4$, $8_{19}$, and $10_{124}$. It is likely that more counterexamples can be identified with a deeper computational sweep.

The inaccessibility of the superhelical configuration in the annealing process highlights the importance of initial conditions. An ensemble of randomized but topologically equivalent initial configurations, for example reached as the result of a Brownian Dynamics simulation of a semiflexible polymer knot, may provide an interesting set of ropelength-volume data and highlight new families of local minima.

The algorithm used in this study was designed to minimize ropelength, not volume. We hope that this work can inspire further investigation into volume minimization of knots. A future algorithm that minimizes volume rather than length would provide a new avenue of investigation in physical knot theory, as the scaling of minimum knot volume with respect to crossing number has not been examined with anywhere near the depth of ropelength data.

\section{Acknowledgements}

We wish to thank Ryan Blair, Hyoungjun Kim and Daniele Celoria for helpful discussions. This material is based upon work supported by the National Science Foundation under Grant No. 2105113.

\bibliographystyle{unsrt}
\bibliography{knotrefshull.bib}

\end{document}